\documentclass[reqno]{amsart}
\usepackage{amsmath}
\usepackage{xspace}
\usepackage{alltt}

\def\eqn#1{(\ref{eq:#1})}
\newcommand{\qBin}[3]{\genfrac{[}{]}{0pt}{0}{#1}{#2}_{#3}}
\numberwithin{equation}{section}
\allowdisplaybreaks[1]

\newcommand{\no}{\noindent}

\newtheorem{theorem}{Theorem}
   \newcommand{\bth}{\begin{theorem}}
   \newcommand{\eth}{\end{theorem}}
\newtheorem{lemma}{Lemma}
   \newcommand{\blem}{\begin{lemma}}
   \newcommand{\elem}{\end{lemma}}

\begin{document}

\makeatother

\title[\tiny{G\"ollnitz-Gordon partitions with weights and parity conditions}]
	{G\"ollnitz-Gordon partitions with weights and parity conditions}

\author[K.~Alladi and A.~Berkovich]
       {Krishnaswami Alladi and Alexander Berkovich}

\address{Department of Mathematics, The University of Florida,
         Gainesville, FL~32611, USA}

\email{alladi@math.ufl.edu}
\email{alexb@math.ufl.edu}


\thanks{
double series representation, infinite hierarchy, Bailey lemma, $q$-trinomial identities}  

\thanks{The first author was supported in part by National Science Foundation Grant DMS--0088975}

\begin{abstract}
A G\"ollnitz-Gordon partition is one in which the parts differ by at least $2$, and 
where the inequality is strict if a part is even.  
Let $Q_i(n)$ denote the number of partitions of $n$ into distinct 
parts $\not\equiv i\pmod{4}$. By attaching weights which are powers 
of $2$ and imposing certain parity conditions on G\"ollnitz-Gordon 
partitions, we show that these are equinumerous with $Q_i(n)$ for 
$i=0,2$. These complement results of G\"ollnitz on $Q_i(n)$ for 
$i=1,3,$ and of Alladi who provided a uniform treatment of all four 
$Q_i(n)$, $i=0,1,2,3,$ in terms of weighted partitions into parts 
differing by $\ge 4$. Our approach here provides a uniform treatment 
of all four $Q_i(n)$ in terms of certain double series representations. 
These double series identities are part of a new infinite hierarchy of 
multiple series identities.  
\end{abstract} 

\maketitle

\section{introduction}
\label{sec:1}

For $i=0,1,2,3,$ let $Q_i(n)$ denote the number of partitions of $n$ 
into distinct parts $\not\equiv i\pmod{4}$. The well known (Little) 
Theorem of G\"ollnitz \cite{6} is:
\bth 
For $i=1,3$, $Q_i(n)$ equals the number of partitions of $n$ into parts 
differing by $\ge2$, where the inequality is strict if a part is odd, 
and the smallest part is $>\frac{(4-i)}{2}$. 
\label{thm:1}
\eth
\no
The analytic representation of Theorem~\ref{thm:1} is
\begin{equation}
\sum^{\infty}_{n=0}\frac{q^{n^2+n}(-q;q^2)_n}{(q^2;q^2)_n}
=(-q^2;q^4)_{\infty}(-q^3;q^4)_{\infty}(-q^4;q^4)_{\infty}
\label{eq:1.1}
\end{equation}
when $i=1$, and 
\begin{equation}
\sum^{\infty}_{n=0}\frac{q^{n^2+n}(-q^{-1};q^2)_n}{(q^2;q^2)_n}
=(-q;q^4)_{\infty}(-q^2;q^4)_{\infty}(-q^4;q^4)_{\infty},
\label{eq:1.2}
\end{equation}
when $i=3$. In \eqn{1.1}, \eqn{1.2}, and in what follows, we have used the standard notation
\begin{equation*}
(a)_n=(a;q)_n=\prod^{n-1}_{j=0}(1-aq^j)
\end{equation*}
for any complex number $a$, and
\begin{equation*}
(a)_\infty=\lim_{n\to\infty}(a)_n=\prod^\infty_{j=0}(1-aq^j),  
\end{equation*}
for $|q|<1$. The products on the right in \eqn{1.1}, \eqn{1.2} are also equal to
\begin{equation*}
\frac{1}{(q^2;q^8)_\infty(q^3;q^8)_\infty(q^7;q^8)_\infty}
\end{equation*}
and
\begin{equation*}
\frac{1}{(q;q^8)_\infty(q^5;q^8)_\infty(q^6;q^8)_\infty},
\end{equation*}
respectively, which have obvious interpretations as generating 
functions of partitions into parts in certain residue classes 
$\pmod{8}$, repetition allowed. 
The equally well known G\"ollnitz-Gordon partition theorem is
\bth 
For $i=1,3$, the number of partitions into parts $\equiv \pm i, 4\pmod{8}$ 
equals the number of partitions into parts differing by $\ge 2$, where 
the inequality is strict if a part is even, and the smallest part is $\ge i$.
\label{thm:2}
\eth
The analytic representation of Theorem~\ref{thm:2} is
\begin{equation}
\sum^{\infty}_{n=0}\frac{q^{n^2}(-q;q^2)_n}{(q^2;q^2)_n}
=\frac{1}{(q;q^8)_{\infty}(q^4;q^8)_{\infty}(q^7;q^8)_{\infty}}
\label{eq:1.3}
\end{equation}
when $i=1$, and 
\begin{equation}
\sum^{\infty}_{n=0}\frac{q^{n^2+2n}(-q;q^2)_n}{(q^2;q^2)_n}
=\frac{1}{(q^3;q^8)_{\infty}(q^4;q^8)_{\infty}(q^5;q^8)_{\infty}}
\label{eq:1.4}
\end{equation}
when $i=3$. Actually \eqn{1.3} and \eqn{1.4} are equations (36) and (37) 
in Slater's famous list \cite{Sl2}, but it was G\"ollnitz \cite{6} and Gordon \cite{7} 
who independently realized their combinatorial interpretation. 

By a reformulation of the (Big) Theorem of G\"ollnitz \cite{6} (not Theorem~\ref{thm:1}) 
using certain quartic transformations, Alladi \cite{1} provided a uniform 
treatment of all four partition functions $Q_i(n),\;i=0,1,2,3$ in terms 
of partitions into parts differing by $\ge 4$, and with certain 
powers of $2$ as weights attached. As a consequence, it was noticed in \cite{1}
that $Q_2(n)$ and $Q_0(n)$ possess certain more interesting properties 
than their well known counterparts $Q_1(n)$ and $Q_3(n)$. In particular, 
$Q_2(n)$ alone among the four functions satisfies the property that 
for every positive integer $k$, $Q_2(n)$ is a multiple of $2^k$ for 
almost all $n$ which was proved by Gordon in an Appendix to \cite{1}. 

Our goal is to prove Theorem~\ref{thm:3} in \S2 which shows that by attaching 
weights which are powers of $2$ to the G\"ollnitz-Gordon partitions 
of $n$, and by imposing certain parity conditions, this is made equal 
to $Q_2(n)$. Here by a {\it G\"ollnitz-Gordon partition} we mean a partition 
into parts differing by $\ge 2$, where the inequality is strict 
if a part is even. There is a similar result for $Q_0(n)$, and this 
is stated as Theorem~\ref{thm:4} at the end of \S2. Theorems~\ref{thm:3} 
and~\ref{thm:4} are nice complements to Theorem~\ref{thm:1} and to results of 
Alladi \cite{1}.

A combinatorial proof of Theorem~\ref{thm:3} is given in full in the next section. 
Theorem~\ref{thm:4} is only stated, and its proof which is similar, is omitted. 

In proving Theorem~\ref{thm:3} we are able to cast it as an analytic identity
(see \eqn{3.2} in \S3) which equates a double series with the product which is 
the generating function of $Q_2(n)$. It turns out that there is a 
two parameter refinement of \eqn{3.2} (see \eqn{3.3} of \S3) which leads to similar 
double series representations for all four products 
\begin{equation*}
\prod_{m>0, m\not\equiv i\pmod{4}}(1+q^m)
\end{equation*}
for $i=0,1,2,3$. It will be shown in \S3 that only in the cases $i=1,3$ 
do these double series reduce to the single series in \eqn{1.1} and \eqn{1.2}. 

Actually, the double series identity \eqn{3.2} is the case $k=2$ of a new 
infinite hierarchy of identities valid for every $k\ge 1$. In \S4 we use 
a limiting case of Bailey's lemma to derive this hierarchy. We give a 
partition theoretic interpretation of the case $k=1$ and state without proof 
a doubly bounded polynomial identity which yields our new hierarchy as 
a limiting case. This polynomial identity will be investigated in detail 
elsewhere.     

\section{a new weighted partition theorem}
\label{sec:2}

Normally, by the parity of an integer we mean its residue class $\pmod{2}$. 
Here by the parity of an odd (or even) integer we mean its residue 
class $\pmod{4}$. 

Next, given a partition $\pi$ into parts differing by $\ge 2$, by a 
{\it{chain}} $\chi$ in $\pi$ we mean a maximal string of parts 
differing by exactly $2$. Thus every partition into parts differing by 
$\ge 2$ can be decomposed into chains. Note that if one part of a 
chain is odd (resp. even), then all parts of the chain are odd 
(resp. even). Hence we may refer to a chain as an odd chain or an 
even chain. Also let $\lambda (\chi)$ denote the least part of 
a chain $\chi$ and $\lambda (\pi)$ the least part of $\pi$. 

Note that in a G\"ollnitz-Gordon partition, since the gap between 
even parts is $>2$, this is the same as saying that every even chain 
is of length $1$, that is, it has only one element. 

Finally, given part $b$ of partition $\pi$, by $t(b;\pi)=t(b)$ we denote 
the number of odd parts of $\pi$ that are $<b$. With this new statistic $t$ 
we now have 
\bth 
Let $\mathcal S$ denote the set of all special G\"ollnitz-Gordon 
partitions, namely, G\"ollnitz-Gordon partitions $\pi$ satisfying 
the parity condition that for every even part $b$ of $\pi$
\begin{equation}
b\equiv 2t(b)\pmod{4}.
\label{eq:2.1}
\end{equation}
Decompose each $\pi \in \mathcal S$ into chains $\chi$ and define the weight 
$\omega (\chi)$ as 
\begin{equation}
\omega (\chi)=
\left\{\begin{array}{l}
2, \mbox{ if $\chi$ is an odd chain, $\lambda (\chi) \ge 5$, 
          and $\lambda (\chi)\equiv 1+2t(\lambda (\chi))\pmod{4}$,} \\
1, \mbox{otherwise}.
\end{array}\right.
\label{eq:2.2}
\end{equation}
The weight $\omega (\pi)$ of the partition $\pi$ is defined multiplicatively 
as 
\begin{equation*}
\omega (\pi)=\prod_{\chi}\omega (\chi),
\end{equation*}
the product over all chains $\chi$ of $\pi$. We then have 
\begin{equation*}
Q_2(n)=\sum_{\pi\in\mathcal S, \, \sigma (\pi)=n}\omega (\pi), 
\end{equation*}
where $\sigma (\pi)$ is the sum of the parts of $\pi$.
\label{thm:3}
\eth

{\bf{Proof:}} Consider the partition $\pi : b_1+b_2+...+b_N$, 
$\pi\in {\mathcal S}$, where contrary to the standard practice of 
writing parts in descending order, we now  have $b_1<b_2<...<b_N.$ 
Subtract $0$ from $b_1$, $2$ from $b_2$, ..., $2N-2$ from $b_N$, to get 
a partition $\pi^*$. We call this process the {\it{Euler subtraction}}. 
Note that in $\pi^*$ the even parts cannot repeat, but the odd parts can. 
Let the parts of $\pi^*$ be $b^*_1\le b^*_2\le ...\le b^*_N.$

Now identify the parts of $\pi$ which are odd, and which are the 
smallest parts of chains and satisfy both the parity and 
low bound conditions in \eqn{2.2}. 
Mark such parts with a {\it{tilde}} at the top. That is, if $b_k$ is 
such a part, we write $b_k=\tilde b_k$ for purposes of identification.
Let $\tilde b_k$ yield $\tilde b^*_k=b^*_k$ after the Euler subtraction. 

Next, split the parts of $\pi^*$ into two piles $\pi^*_1$ and 
$\pi^*_2$, with $\pi^*_1$ consisting only of certain odd parts, 
and $\pi^*_2$ containing the remaining parts. In this decomposition 
we adopt the following rule:

(a) the odd parts of $\pi^*$ which are not identified as above are 
put in $\pi^*_1$. 

(b) the odd parts of $\pi^*$ which have been identified could 
be put in either $\pi^*_1$ or $\pi^*_2$. 

Thus we have two choices for each identified part. 

Let us say, in a certain given situation, after making the choices, 
we have $n_1$ parts in $\pi^*_1$ and $n_2$ parts in $\pi^*_2$. We 
now add $0$ to the smallest part of $\pi^*_2$, $2$ to the second smallest 
part of $\pi^*_2$, ..., $2n_2-2$ to the largest part of $\pi^*_2$, 
$2n_2$ to the smallest part of $\pi^*_1$, $2n_2+2$ to the second 
smallest part of $\pi^*_1$, ..., $2(n_1+n_2)-2=2N-2$ to the largest 
part of $\pi^*_1$. We call this the {\it{Bressoud redistribution}} process. 
As a consequence of this redistribution, we have created two partitions 
$\pi_1$ (out of $\pi^*_1$) and $\pi_2$ (out of $\pi^*_2$) satisfying 
the following conditions:

(i) $\pi_1$ consists only of distinct odd parts, with each odd part being 
greater than twice the number of parts of $\pi_2$. 

(ii) Since both the even and odd parts of $\pi^*_2$ are distinct, the 
parts of $\pi_2$ differ by $\ge 4$. Also since the odd parts of $\pi^*_2$ 
are chosen from the smallest of parts of certain chains in $\pi$, the 
odd parts of $\pi_2$ actually differ by $\ge 6$, and each such odd part 
is $\ge 5$.

In transforming the original partition $\pi$ into the pair $(\pi_1,\pi_2)$, 
we need to see how the parity conditions of $\pi$ given by \eqn{2.1} and 
\eqn{2.2} transform to parity conditions in $\pi_1$ and $\pi_2$. 

First observe that since the parity conditions on $\pi$ are imposed 
only on the even parts of $\pi$ and the {\it{identified}} odd parts 
of $\pi$, the transformed parity conditions (to be determined below) 
will be imposed only on $\pi_2$ and not on $\pi_1$. Thus $\pi_1$ will 
satisfy only condition (i) above. 

Suppose $b_k$ is an even part of $\pi$ and that $t(b_k;\pi)=t$, that is 
there are $t$ odd parts of $\pi$ which are less than $b_k$. Now $b_k$  
becomes 
\begin{equation*}
b^*_k=b_k-(2k-2)
\end{equation*}
after the Euler subtraction. Notice that $t(b^*_k;\pi^*)=t(b_k;\pi)=t$. 
Now suppose that from among the $t$ odd parts of $\pi^*$ less than 
$b^*_k$, $r$ of them are put in $\pi^*_1$ and the remaining $t-r$ odd 
parts are put in $\pi^*_2$. Then $b^*_k$ becomes the $(k-r)-th$ smallest 
part in $\pi^*_2$. So in the Bressoud redistribution process, 
$2(k-r)-2$ is added to $b^*_k$ making it a new even part $e_{k-r}$ in 
$\pi_2$. Thus 
\begin{equation}
e_{k-r}=b^*_k+2(k-r)-2=b_k-(2k-2)+2(k-r)-2=b_k-2r.
\label{eq:2.3}
\end{equation}  
We see from \eqn{2.1} and \eqn{2.3} that 
\begin{equation}
e_{k-r}\equiv 2t-2r=2(t-r)=2t(e_{k-r};\pi_2)\pmod{4}
\label{eq:2.4}
\end{equation}
and so the parity condition \eqn{2.1} on the even parts does not 
change when going to $\pi_2$. Thus we may write \eqn{2.4} in short as 
\begin{equation}
e\equiv 2t(e)\pmod{4}
\label{eq:2.5}
\end{equation}
for any even part in $\pi_2$. 

Now we need to determine the parity conditions on the odd parts in 
$\pi_2$ which are derived from some of the identified odd parts of $\pi$. 
To this end suppose that $\tilde b_k$ is an identified odd part of $\pi$ 
which becomes $\tilde b^*_k=\tilde b_k-(2k-2)$ in $\pi^*$ due to the 
Euler subtraction, and that $\tilde b^*_k$ is placed in $\pi^*_2$. 
Let $t(\tilde b_k;\pi)=t$. Notice that 
\begin{equation*}
t(\tilde b_k;\pi)=t(\tilde b^*_k;\pi^*)=t.
\end{equation*}
Suppose that from among the $t$ odd parts of $\pi^*$ which are 
$\tilde b^*_k$, $r$ of them are placed in $\pi^*_1$ and the remaining 
$t-r$ are placed in $\pi^*_2$. Then $\tilde b^*_k$ becomes the 
$(k-r)-th$ smallest part in $\pi^*_2$. Thus under the Bressoud 
redistribution, $2(k-r)-2$ is added to it to yield the part $f_k$ given by
\begin{equation*}
f_k=\tilde b^*_k+2(k-r)-2=\tilde b_k-(2k-2)+(2(k-r)-2)=\tilde b_k-2r
\end{equation*}
as in \eqn{2.3}. Therefore the parity condition \eqn{2.2} yields
\begin{equation*}
f_k\equiv 1+2t-2r=1+2(t-r)\pmod{4}.
\end{equation*}
But $t(f_k;\pi_2)=t-r$. So this could be expressed in short as 
\begin{equation}
f\equiv 1+2t(f)\pmod{4}
\label{eq:2.6}
\end{equation}
for any odd part of $\pi_2$. Thus the pair of partitions $(\pi_1,\pi_2)$ 
is determined by condition (i) on $\pi_1$, and conditions (ii) and 
the parity conditions \eqn{2.5} and \eqn{2.6} on $\pi_2$. 

In going from $\pi$ to the pair $(\pi_1,\pi_2)$ we had a {\it{choice}} 
of deciding whether an {\it{identified}} part of $\pi$ would end up 
in $\pi_1$ or $\pi_2$. This choice is precisely the weight 
$\omega (\chi)=2$ associated with certain chains $\chi$. The weight 
of the partition $\pi$ is computed mutiplicatively because these choices 
are {\it{independent}}. So what we have established up to now is:

\blem
The weighted count of the special G\"ollnitz-Gordon partitions of $n$ 
equals the number of bipartitions $(\pi_1,\pi_2)$ of $n$ satisfying 
conditions (i), (ii), \eqn{2.5} and \eqn{2.6}.
\label{lem:1}
\elem

Next, we discuss a bijective map 
\begin{equation}
\pi_2 \mapsto (\pi_3,\pi_4),
\label{eq:2.7}
\end{equation}
where $\pi_3$ is a partition into distinct multiples of $4$ and $\pi_4$ 
is a partition into distinct odd parts such that 
\begin{equation}
\nu(\pi_2)=\nu(\pi_3)
\label{eq:2.8}
\end{equation}
and
\begin{equation}
2\nu(\pi_2)>\Lambda(\pi_4).
\label{eq:2.9}
\end{equation}
Here by $\nu(\pi)$ we mean the number of parts of a partition $\pi$ and 
by $\Lambda(\pi)$ the largest part of $\pi$. 

To describe the map \eqn{2.7} we represent $\pi_2$ as a Ferrers graph
with weights $1,2$ or $4$, at each node. We construct the graph as 
follows: 

1) With each odd (resp. even) part f (resp. e) of $\pi_2$
   we associate a row of $\frac{3+f+2t(f)}{4}$ (resp. $\frac{e+2t(e)}{4}$) nodes.

2) We place a $1$ at end of any row that represents an odd part of $\pi_2$. 

3) Every node in the column directly above each $1$ is given weight $2$. 

4) Each remaining node is given weight $4$. 

Every part of $\pi_2$ is given by the sum of weights in an associated row.
It is clear from these weights, that the partition represented by 
this weighted Ferrers graph  satisfies precisely the conditions 
(ii), \eqn{2.5} and \eqn{2.6} that characterize $\pi_2$.      

$$
\begin{array}{lllllllllll}
\pi_2: & & 4 & 2 & 4 & 4 & 2 & 4 & 4 & 4 & 1 \\
       & & 4 & 2 & 4 & 4 & 2 & 4 & 4 &   &   \\
       & & 4 & 2 & 4 & 4 & 1 &   &   &   &   \\
       & & 4 & 1 &   &   &   &   &   &   &
\end{array}
$$

Next we extract from this weighted Ferrers graph all columns with 
a $1$ at the bottom, and assemble these columns as rows to form a $2$-modular 
Ferrers graph as shown below. 

$$
\begin{array}{llllll}
\pi_4: & & 2 & 2 & 2 & 1 \\
       & & 2 & 2 & 1 &   \\
       & & 1 &   &   &
\end{array}
$$

Clearly this $2$-modular graph represents a partition $\pi_4$ that 
satisfies condition \eqn{2.9}. 

After this extraction, the decorated graph of $\pi_2$ becomes a $4$-modular 
graph (in this case a graph with weight $4$ at every node). This graph 
$\pi_3$ clearly satisfies \eqn{2.8}. 

$$
\begin{array}{llllllll}
\pi_3: & & 4 & 4 & 4 & 4 & 4 & 4 \\
       & & 4 & 4 & 4 & 4 & 4 &   \\
       & & 4 & 4 & 4 &   &   &   \\
       & & 4 &   &   &   &   &
\end{array}
$$

It is easy to check that \eqn{2.7} is a bijection. Thus Lemma~\ref{lem:1} can be recasted 
in the form 
\blem
The weighted count of the special G\"ollnitz-Gordon partitions of $n$ 
as in Theorem~\ref{thm:3} is equal to the number of partitions of $n$ in the form 
$(\pi_1, \pi_3, \pi_4)$ where

(iii) $\pi_3$ consists only of distinct multiples of 4,

(iv) $\pi_4$ has distinct odd parts and $\Lambda(\pi_4)<2\nu(\pi_3)$, 

(v) $\pi_1$ has distinct odd parts and $\lambda(\pi_1)>2\nu(\pi_3)$,
\label{lem:2}   
\elem
 
Finally, observe that conditions (iv) and (v) above yield partitions 
into distinct odd parts (without any other conditions). This together 
with (iii) yields partitions counted by $Q_2(n)$, thereby completing the 
combinatorial proof of Theorem~\ref{thm:3}. 

In a similar fashion, we can obtain the following representation for 
$Q_0(n)$ with weights and parity conditions imposed on the 
G\"ollnitz-Gordon partitions:
\bth 
Let ${\mathcal S}^*$ denote the set of all special G\"ollnitz-Gordon 
partitions, namely, G\"ollnitz-Gordon partitions $\pi$ satisfying 
the parity condition that for every even part $b$ of $\pi$
\begin{equation}
b\equiv 2(t(b)-1)\pmod{4}.
\label{eq:2.10}   
\end{equation}
Decompose each $\pi \in {\mathcal S}^*$ into chains $\chi$ and define the weight 
$\omega (\chi)$ as 
\begin{equation}
\omega (\chi)=
\left\{\begin{array}{l}
2, \mbox{ if $\chi$ is an odd chain, $\lambda (\chi) \ge 3$, 
          and $\lambda (\chi)\equiv 2t(\lambda (\chi))-1\pmod{4}$,} \\
1, \mbox{otherwise}.
\end{array}\right.
\label{eq:2.11}   
\end{equation}
The weight $\omega (\pi)$ of the partition $\pi$ is defined multiplicatively 
as 
\begin{equation*}
\omega (\pi)=\prod_{\chi}\omega (\chi),
\end{equation*}
the product over all chains $\chi$ of $\pi$. We then have 
\begin{equation*}
Q_0(n)=\sum_{\pi \in {\mathcal S}^*, \, \sigma (\pi)=n}\omega (\pi), 
\end{equation*}
where $\sigma (\pi)$ is the sum of the parts of $\pi$.
\label{thm:4}   
\eth

\section{series representations}
\label{sec:3}

If we let $\nu(\pi_1)=n_1$ and $\nu(\pi_2)=n_2$, then \eqn{2.7} and 
conditions (iii), (iv), and (v) of Lemma~\ref{lem:2} imply that the generating 
function of all such triples of partitions $(\pi_1, \pi_3, \pi_4)$ is 
\begin{equation}
\frac{q^{n^2_1+2n_1n_2}}{(q^2;q^2)_{n_1}}.
\frac{q^{2n^2_2+2n_2}}{(q^4;q^4)_{n_2}}.(-q;q^2)_{n_2}.
\label{eq:3.1}   
\end{equation} 
If the expression in \eqn{3.1} is summed over all non-negative integers $n_1$ 
and $n_2$, it yields
\begin{equation*}
\sum_{n_1}\sum_{n_2}\frac{q^{n^2_1+2n_1n_2+2n^2_2+2n_2}(-q;q^2)_{n_2}}
{(q^2;q^2)_{n_1}(q^4;q^4)_{n_2}}
=\sum_{n_2}\frac{q^{2n^2_2+2n_2}(-q;q^2)_{n_2}}{(q^4;q^4)_{n_2}}
\sum_{n_1}\frac{q^{n^2_1+2n_1n_2}}{(q^2;q^2)_{n_1}}
\end{equation*}
\begin{equation*}
=\sum_{n_2}\frac{q^{2n^2_2+2n_2}(-q;q^2)_{n_2}}{(q^4;q^4)_{n_2}}
(-q^{2n_2+1};q^2)_\infty
=(-q;q^2)_\infty\sum_{n_2}\frac{q^{2n^2_2+2n_2}}{(q^4;q^4)_{n_2}}
\end{equation*}
\begin{equation}
=(-q;q^2)_\infty(-q^4;q^4)_\infty
=(-q;q^4)_\infty(-q^3;q^4)_\infty(-q^4;q^4)_\infty
=\sum_nQ_2(n)q^n.
\label{eq:3.2}   
\end{equation}
By just following the above steps we can actually get a two parameter 
refinement of \eqn{3.2}, namely,
\begin{equation} 
\sum_{n_1}\sum_{n_2}
\frac{z^{n_1}{\omega}^{n_2}q^{n^2_1+2n_1n_2+2n^2_2+2n_2}(-zq;q^2)_{n_2}}
{(q^2;q^2)_{n_1}(q^4;q^4)_{n_2}}
=(-zq;q^4)_\infty(-zq^3;q^4)_\infty(-\omega q^4;q^4)_\infty.
\label{eq:3.3}   
\end{equation}

One may view \eqn{3.2} as the analytic version of Theorem~\ref{thm:3}. In 
reality, the correct way to view \eqn{3.2} is that, if the summand on the 
left is decomposed into three factors as \eqn{3.1}, then \eqn{3.2} is the analytic 
version of the statement that the number of partitions of an integer 
$n$ into the triple of partitions $(\pi_1, \pi_3, \pi_4)$ is equal to 
$Q_2(n)$. This is of course only the final step of the proof given above. 
and \eqn{3.2}, which is quite simple, is equivalent to it.  

The advantage in the two parameter refinement \eqn{3.3} is that by suitable 
choice of the parameters we get similar representations involving $Q_i(n)$ 
for $i=0,1,3$. For example, if we replace $\omega$ by $\omega q^{-2}$ in 
\eqn{3.3} we get
\begin{equation} 
\sum_{n_1}\sum_{n_2}
\frac{z^{n_1}{\omega}^{n_2}q^{n^2_1+2n_1n_2+2n^2_2}(-zq;q^2)_{n_2}}
{(q^2;q^2)_{n_1}(q^4;q^4)_{n_2}}
=(-zq;q^4)_\infty(-zq^3;q^4)_\infty(-\omega q^2;q^4)_\infty,
\label{eq:3.4}   
\end{equation}
which is the analytic representation of Theorem~\ref{thm:4} above.

Next, replacing $z$ by $zq$ and $\omega$ by $\omega q^{-1}$ in \eqn{3.3} we get
\begin{equation} 
\sum_{n_1}\sum_{n_2}
\frac{z^{n_1}{\omega}^{n_2}q^{n^2_1+2n_1n_2+2n^2_2+n_1+n_2}(-zq^2;q^2)_{n_2}}
{(q^2;q^2)_{n_1}(q^4;q^4)_{n_2}}
=(-zq^2;q^4)_\infty(-\omega q^3;q^4)_\infty(-zq^4;q^4)_\infty.
\label{eq:3.5}   
\end{equation}
Now choose $z=1$ in \eqn{3.5}. Then the double series on the left becomes 
\begin{equation}
\sum_{n_1}\sum_{n_2}
\frac{{\omega}^{n_2}q^{n^2_1+2n_1n_2+2n^2_2+n_1+n_2}}
{(q^2;q^2)_{n_1}(q^2;q^2)_{n_2}}
=\sum_{n_1}\sum_{n_2}
\frac{{\omega}^{n_2}q^{(n_1+n_2)^2+n^2_2+n_1+n_2}}
{(q^2;q^2)_{n_1}(q^2;q^2)_{n_2}}.
\label{eq:3.6}   
\end{equation}
If we now put $n=n_1+n_2$ and $j=n_2$, then \eqn{3.6} could be rewritten in 
the form
\begin{equation*}
\sum_n\frac{q^{n^2+n}}{(q^2;q^2)_n}
\sum^n_{j=0}\frac{{\omega}^jq^{j^2}(q^2;q^2)_n}{(q^2;q^2)_j(q^2;q^2)_{n-j}}
\end{equation*}
\begin{equation}
=\sum^{\infty}_{n=0}\frac{q^{n^2+n}(-\omega q;q^2)_n}{(q^2;q^2)_n}
=(-q^2;q^4)_{\infty}(-\omega q^3;q^4)_{\infty}(-q^4;q^4)_{\infty},
\label{eq:3.7}   
\end{equation}
which is the single series identity \eqn{1.1} in a refined form.

Similarly, replacing $\omega$ by $\omega q^{-3}$ and $z$ by $zq$ in \eqn{3.2} 
we get
\begin{equation} 
\sum_{n_1}\sum_{n_2}
\frac{z^{n_1}{\omega}^{n_2}q^{n^2_1+2n_1n_2+2n^2_2+n_1-n_2}(-zq^2;q^2)_{n_2}}
{(q^2;q^2)_{n_1}(q^4;q^4)_{n_2}}
=(-zq^2;q^4)_\infty(-\omega q;q^4)_\infty(-zq^4;q^4)_\infty.
\label{eq:3.8}   
\end{equation}
Now the choice $z=1$ makes the double series in \eqn{3.8} as
\begin{equation}
\sum_{n_1}\sum_{n_2}
\frac{{\omega}^{n_2}q^{n^2_1+2n_1n_2+2n^2_2+n_1-n_2}}
{(q^2;q^2)_{n_1}(q^2;q^2)_{n_2}}
=\sum_{n_1}\sum_{n_2}
\frac{{\omega}^{n_2}q^{(n_1+n_2)^2+n^2_2+n_1-n_2}}
{(q^2;q^2)_{n_1}(q^2;q^2)_{n_2}}.
\label{eq:3.9}   
\end{equation}
Once again, putting $n=n_1+n_2$ and $j=n_2$ makes \eqn{3.9} into  
\begin{equation*}
\sum_n\frac{q^{n^2+n}}{(q^2;q^2)_n}
\sum^n_{j=0}\frac{{\omega}^jq^{j^2-2j}(q^2;q^2)_n}{(q^2;q^2)_j(q^2;q^2)_{n-j}}
\end{equation*}
\begin{equation}
=\sum^{\infty}_{n=0}\frac{q^{n^2+n}(-\omega q^{-1};q^2)_n}{(q^2;q^2)_n}
=(-q^2;q^4)_{\infty}(-\omega q;q^4)_{\infty}(-q^4;q^4)_{\infty},
\label{eq:3.10}   
\end{equation}
which is a refinement of the single series identity \eqn{1.2}. Thus 
precisely in the cases $i=1,3$, can the double series be reduced to 
single series by setting one of the parameters $z=1$. 

\section{a new infinite hierarchy}
\label{sec:4}

Identity \eqn{3.2} given above is just the case $k=2$ of a new infinite 
hierarchy of multiple series identities \eqn{4.12} given below.

To derive this hierarchy, we will need the definition of a Bailey pair, 
and a special case of Bailey's lemma which produces a new Bailey pair 
from a given Bailey pair \cite{A1}. 

{\bf{Definition}}: A pair of sequences $\alpha_n(q), \beta_n(q)$ is called 
a Bailey pair (relative to $1$) if for all $n\ge 0$ 
\begin{equation}
\beta_n(q)=\sum^n_{i=0} \frac{\alpha_i(q)}{(q)_{n-i}(q)_{n+i}}.
\label{eq:4.1} 
\end{equation}
By setting $a=1, \rho_1=-q^{\frac{1}{2}}$, and letting $\rho_2\to\infty$ in the 
formulas (3.29) and (3.30) of \cite{A1}, we obtain the following limiting 
case of Bailey's lemma:
\blem
Suppose $(\alpha_n(q), \beta_n(q))$ is a Bailey pair. 
Then $(\alpha_n^{(1)}(q), \beta_n^{(1)}(q))$ is another Bailey pair, where  
\begin{equation}
\alpha_n^{(1)}(q)=q^{\frac{n^2}{2}}\alpha_n(q),
\label{eq:4.2} 
\end{equation}
\begin{equation}
\beta_n^{(1)}(q)=
\sum^n_{i=0} \frac{(-\sqrt q)_i}{(q)_{n-i}(-\sqrt q)_n}q^{\frac{i^2}{2}}\beta_i(q).
\label{eq:4.3} 
\end{equation}
\label{lem:3} 
\elem
From $(\alpha_n^{(1)}(q),\beta_n^{(1)}(q))$ one can produce next Bailey pair 
$(\alpha_n^{(2)}(q),\beta_n^{(2)}(q))$ simply using 
$(\alpha_n^{(1)}(q),\beta_n^{(1)}(q))$ as the initial Bailey pair. 
It is easy to check that the $k$-fold iteration of (the limiting case of) Bailey's Lemma yields
\begin{equation}
\alpha_n^{(k)}(q)=q^{k\frac{n^2}{2}}\alpha_n(q),
\label{eq:4.4} 
\end{equation}
\begin{equation}
\beta_n^{(k)}(q)=
\sum_{\overrightarrow n} \frac{q^{\frac{N_1^2+N_2^2+\cdots +N_k^2}{2}}(-\sqrt q)_{n_k}}
{(q)_{n-N_1} (-\sqrt q)_n (q)_{n_1} \cdots (q)_{n_{k-1}}} \beta_{n_k}(q),
\label{eq:4.5} 
\end{equation}
where $\overrightarrow n=(n_1,n_2,\cdots,n_k)$ and $N_i=n_i+n_{i+1}+\cdots +n_k,$ 
with $i=1,2,\ldots,k$.
In \cite{Sl1}, \cite{Sl2} Slater derived A-M families of Bailey pairs to produce the celebrated 
list of $130$ identities of the Rogers-Ramanujan type. We shall need her $E(4)$ pair:
\begin{eqnarray}
\alpha_n & = & 
\left\{\begin{array}{l}
(-1)^n q^{n^2}(q^n+q^{-n}), \mbox{if }n>0, \notag \\
1, \mbox{if } n=0,
\end{array}\right.\notag\\
\beta_n & = & \frac{q^n}{(q^2;q^2)_n}.
\label{eq:4.6} 
\end{eqnarray}

\no
It follows from \eqn{4.1} and \eqn{4.4} - \eqn{4.6} that
\begin{equation}
\sum_{\overrightarrow n} \frac{q^{\frac{1}{2}(N_1^2+N_2^2+\cdots +N_k^2)+N_k}(-\sqrt q)_{n_k}}
{(q)_{n-N_1}(q)_{n_1}(q)_{n_2} \cdots (q)_{n_{k-1}}(q^2;q^2)_{n_k}}= 
\frac{(-\sqrt q)_n}{(q)_n}\sum_{j=-n}^n(-1)^jq^{\frac{k+2}{2}j^2+j}
\qBin{2n}{n+j}{q},
\label{eq:4.7} 
\end{equation}
where $q$-binomial coefficients are defined as 
\begin{equation}
\qBin{n+m}{n}{q}=\frac{(q^{m+1})_n}{(q)_n}.
\label{eq:4.8} 
\end{equation}
It is easy to check that
\begin{equation}
\lim_{n\to\infty}\qBin{n}{m}{q}=\frac{1}{(q)_m},
\label{eq:4.9} 
\end{equation}
and
\begin{equation}
\lim_{n\to\infty} \qBin{2n}{n+j}{q}=\frac{1}{(q)_\infty}.
\label{eq:4.10} 
\end{equation}

Next, we recall Jacobi's triple product identity
\begin{equation}
\sum_{n=-\infty}^\infty q^{n^2}z^n=(q^2,-qz,-\frac{q}{z};q^2)_\infty,
\label{eq:4.11} 
\end{equation}
where $(a_1,a_2,\ldots,a_m;q)_\infty=(a_1)_\infty(a_2)_\infty\ldots(a_m)_\infty$.\\
If we let $n$ tend to infinity in \eqn{4.7} with $q\to q^2$, we obtain with the aid
of \eqn{4.10} and \eqn{4.11} the desired identity
\begin{eqnarray}
& \sum_{\overrightarrow n} & \frac{q^{N_1^2+\cdots +N_k^2+2N_k}(-q;q^2)_{n_k}}
{(q^2;q^2)_{n_1}\ldots(q^2;q^2)_{n_{k-1}}(q^4;q^4)_{n_k}} \nonumber \\
& = & \frac{(-q;q^2)_\infty}{(q^2;q^2)_\infty}(q^{2k+4},q^k,q^{k+4};q^{2k+4})_\infty\nonumber\\
& = & \frac{(q^2;q^4)_\infty}{(q)_\infty}(q^{2k+4},q^k,q^{k+4};q^{2k+4})_\infty. 
\label{eq:4.12} 
\end{eqnarray}

\no
Here we used the simple relation
\begin{equation*}
\frac{(q^2;q^4)_\infty}{(q)_\infty}=\frac{(q,-q;q^2)_\infty}{(q,q^2;q^2)_\infty}=
\frac{(-q;q^2)_\infty}{(q^2;q^2)_\infty}.
\end{equation*}
Making use of
\begin{equation}
\frac{(-q;q^2)_\infty}{(q^2;q^2)_\infty}(q^8,q^2,q^6;q^8)=
\frac{(-q;q^2)_\infty}{(q^4;q^8)_\infty}=
(-q;q^2)_\infty(-q^4;q^4)_\infty,
\label{eq:4.13}
\end{equation}
it is straightforward to verify that \eqn{4.12} with $k=2$ yields \eqn{3.2}, as claimed.

\no
When $k=1$, \eqn{4.12} becomes
\begin{equation}
\sum_{n\ge 0}\frac{q^{n^2+2n}(-q;q^2)_n}{(q^4;q^4)_n}=
\frac{(q^2;q^4)_\infty(q^6,q^1,q^5;q^6)_\infty}{(q)_\infty}=
\frac{(-q^3;q^6)_\infty}{(q^4,q^8;q^{12})_\infty}. 
\label{eq:4.14} 
\end{equation}
Surprisingly, \eqn{4.14} is missing from the Slater list. It was given by Andrews in \cite{A2}.

By using the statistic $s(b;\pi)$ = number of even parts of the
partition $\pi$ which are less than the part $b$, it can be shown that 
the the following partition theorem is a combinatorial interpretation 
of (4.14):

\bth
Let $G(N)$ denote the number of partitions $\pi$ of 
$N$ into distinct parts such that no gap between consecutive parts is 
$\equiv 1 \pmod{4}$,  and
where the $k$-th smallest part $b$ is $\equiv 1+2k+2s(b;\pi)\pmod{4}$ 
if  $b$ is odd, and $\equiv 2+2k+2s(b;\pi) \pmod{4}$, if  $b$  is even. \\
Let $P(N)$ denote the number of partitions of $N$ into parts $\equiv\pm 3,\pm 4 \pmod{12}$, 
such that parts $\equiv 3 \pmod{6}$ are distinct. Then,
\begin{equation*}
G(N)=P(N).
\end{equation*}
\label{thm:5} 
\eth
{\it{Remark}}: Theorem~\ref{thm:5} can be stated without appeal to the statistic 
$s(b;\pi)$, but we preferred to state it this way to emphasise a
different parity condition and to show similarity with Theorems~\ref{thm:3} and  ~\ref{thm:4}. 

It would be interesting to find partition theoretical interpretation of \eqn{4.12} with $k>2$. 
To this end we observe that the product on the right of \eqn{4.12} with $k\equiv 0 \pmod{4}$ 
can be interpreted as a generating function for partitions into 
parts $\not\equiv 2\pmod{4}$, $\not\equiv 0,\pm k \pmod{2k+4}$. 

\no
It is instructive to compare this product
\begin{equation*}
\prod_{\substack{n\ge 1 \\ n\not\equiv 2\pmod{4}\\ n\not\equiv 0,\pm(2K-2)\pmod{4K}}}
(1-q^n)^{-1}
\end{equation*}
and the generalized G\"ollnitz-Gordon product ((7.4.4); \cite{4})
\begin{equation*}
\prod_{\substack{n\ge 1\\n\not\equiv 2\pmod{4}\\n\not\equiv 0,\pm(2\tilde k-1)\pmod{4\tilde k}}}
(1-q^n)^{-1}.
\end{equation*}
Here $K=1+\frac{k}{2}$ with $k\equiv 0\pmod{4}$ and $\tilde k$ is a positive integer.

\no
The right hand side of \eqn{4.12} can be rewritten as 
\begin{equation*}
\frac{(q^2;q^4)_\infty(q^k,q^{k+4},q^{k+2},-q^{k+2};q^{2k+4})_\infty
(q^{4k+8};q^{4k+8})_\infty}{(q)_\infty},
\end{equation*}
if $k$ is odd, and
\begin{equation*}
\frac{(q^2;q^4)_\infty(q^{2k+4};q^{2k+4})_\infty
(q^{\frac{k}{2}},-q^{\frac{k}{2}},q^{2+\frac{k}{2}},-q^{2+\frac{k}{2}};q^{k+2})_\infty}
{(q)_\infty},
\end{equation*}
if $k\equiv 2\pmod{4}$.

\no 
This enables us to interprete the right hand side of \eqn{4.12} as:\\
A. $k\equiv 1\pmod{2}$. RHS \eqn{4.12} is the generating function for partitions into 
parts $\not\equiv 2\pmod{4}$, $\not\equiv\pm k\pmod{2k+4}$, $\not\equiv 0\pmod{4k+8}$, 
such that parts $\equiv k+2\pmod{2k+4}$ are distinct.
 
\no
B. $k\equiv 2\pmod{4}$. RHS \eqn{4.12} is the generating function for partitions into 
parts $\not\equiv 2\pmod{4}$, $\not\equiv 0\pmod{2k+4}$, such that parts 
$\not\equiv\pm\frac{k}{2}\pmod{k+2}$ are distinct. 

We would like to conclude with the following observation. The hierarchy \eqn{4.12} 
follows in the limit $l,m\to\infty$ from the doubly bounded polynomial identity
\begin{equation}
\begin{split} 
&\sum_{\overrightarrow n,s} q^{N_1^2+\cdots +N_k^2+s^2+2N_k}
\qBin{l+m-N_1}{m-N_1}{q^2}
\prod_{j=1}^{k-1}\qBin{l-\sum_{i=1}^j N_i+n_j}{n_j}{q^2}
\qBin{n_k+\lfloor\frac{l-1-\sum_{i=1}^k N_i-s}{2}\rfloor}{n_k}{q^4}
\qBin{n_k}{s}{q^2}=\\
&\sum_{j=-\infty}^\infty \left\{q^{(4k+8)j^2+4j}\tilde U(l,m,2(k+2)j+1,2j,q^2)-
q^{(4k+8)j^2+4(k+1)j+k}\tilde U(l,m,2(k+2)j+k+1,2j+1,q^2)\right\},
\end{split}
\label{eq:4.15}
\end{equation} 
where $\lfloor x\rfloor$ is the largest integer $\leq x$, $\tilde U(l,m,a,b,q)=
T_w(l,m,a,b,q)+T_w(l,m,a+1,b,q)$, and the refined 
$q$-trinomial coefficients \cite{W} are defined as
\begin{equation}
T_w(l,m,a,b,q):=\sum_{\substack{n=0 \\ n+l\equiv a\pmod{2}}}^l q^{\frac{n^2}{2}}
\qBin{m}{n}{q}\qBin{m+b+\frac{l-a-n}{2}}{m+b}{q}\qBin{m-b+\frac{l+a-n}{2}}{m-b}{q}.
\label{eq:4.16}
\end{equation} 
Using \eqn{4.9} together with Warnaar's limiting formula ((2.26); \cite{W})
\begin{equation}
\lim_{l\to\infty}\tilde U(l,m,a,b,q)=\frac{(-\sqrt q)_m}{(q)_{2m}}\qBin{2m}{m+b}{q},
\label{eq:4.17}
\end{equation}
we obtain \eqn{4.7} with $q\to q^2$ and $n\to m$ as $l\to\infty$ in \eqn{4.12}.
On the other hand, if we let $m\to\infty$ in \eqn{4.16} we find that
\begin{equation}
\lim_{m\to\infty}T_w(l,m,a,b,q)=\frac{1}{(q)_l}T_{AB}(l,a,q),
\label{eq:4.18}
\end{equation}
where the Andrews-Baxter $q$-trinomial coefficients \cite{AB} are defined as
\begin{equation}
T_{AB}(l,a,q):=\sum_{\substack{n=0 \\ n+l\equiv a\pmod{2}}}^l q^{\frac{n^2}{2}}
\qBin{l}{n}{q}\qBin{l-n}{\frac{l-a-n}{2}}{q}.
\label{eq:4.19}
\end{equation}
And so, \eqn{4.15} becomes in the limit $m\to\infty$
\begin{equation}
\begin{split} 
&\sum_{\overrightarrow n,s} q^{N_1^2+\cdots +N_k^2+s^2+2N_k}
\prod_{j=1}^{k-1}\qBin{l-\sum_{i=1}^j N_i+n_j}{n_j}{q^2}
\qBin{n_k+\lfloor\frac{l-1-\sum_{i=1}^k N_i-s}{2}\rfloor}{n_k}{q^4}
\qBin{n_k}{s}{q^2}=\\
&\sum_{j=-\infty}^\infty \left\{q^{(4k+8)j^2+4j}U(l,2(k+2)j+1,q^2)-
q^{(4k+8)j^2+4(k+1)j+k}U(l,2(k+2)j+k+1,q^2)\right\},
\end{split}
\label{eq:4.20}
\end{equation} 
where 
\begin{equation}
U(l,a,q)=T_{AB}(l,a,q)+T_{AB}(l,a+1,q).
\label{eq:4.21}
\end{equation} 
The proof of \eqn{4.15} will be given elsewhere.

{\bf{Acknowledgement}}: We would like to thank Frank Garvan for many stimulating 
discussions and for his help with the diagrams.

\end{document}